\documentclass[12pt]{article}
\usepackage{amssymb}
\usepackage{amsmath}
\usepackage{graphicx}
\newcommand{\R}{{\mathbb R}}
\newcommand{\la}{\lambda}
\newcommand{\lap}{\mbox{$\bigtriangleup$}}

\newcommand{\ra}{{\mbox{$\rightarrow$}}}
\newcommand{\be}{\begin{equation}}
\newcommand{\ee}{\end{equation}}
\newtheorem{mrem}{Remark}
\newtheorem{mthm}{Theorem}

\newtheorem{lem}{Lemma}[section]
\newtheorem{cor}{Corollary}[section]

\begin{document}
\title{A Hopf's Lemma and the Boundary Regularity for the Fractional P-Laplacian}

\author{  Lingyu Jin \thanks{Partially supported by .}
\quad  Yan Li \thanks{ Corresponding author.}   }

\date{\today}
\maketitle
\begin{abstract} We begin the paper with a Hopf's lemma
for a fractional p-Laplacian problem on a half-space. Specifically speaking, we show that the derivative  of the solution
 along the outward normal vector is strictly positive
on the boundary of the half-space. Next we show that positive solutions to a fractional p-Laplacian equation possess certain  H\"{o}lder continuity  up to the boundary.

\end{abstract}
\bigskip

{\bf Key words} Fractional p-Laplacian; Dirichlet problem; Hopf's lemma; Boundary regularity
\bigskip

\section{Introduction}

The strong maximum principle of Eberhard Hopf, often known as the Hopf's lemma\cite{2}, is a classical and  fundamental result
of the theory of second order elliptic partial differential equations.
Its main idea is that \emph{if a function satisfies a second order partial differential inequality of a certain kind in a domain of $R^n$ and attains a maximum in the domain then the function is constant}.
The Hopf's lemma has been generalized to describe the behavior of the solution to an elliptic problem as it approaches a point on the boundary where its maximum is attained.

Its history can be first
traced back to the maximum principle for harmonic functions. In the past decade this lemma has been generalized
as  the  strong maximum principle for singular quasi-linear elliptic differential inequalities(\cite{7}). For a while it was thought that the Hopf's maximum principle applies only to linear differential operators.
 In the later sections of his original paper, however, Hopf considered a more general situation which permits certain nonlinear operators and, in some cases, leads to uniqueness statements in the Dirichlet problem for the mean curvature operator and the Monge-Amp\`{e}re equation.

In the first part of this paper, we prove a Hopf's lemma for a nonlinear non-local pseudo-differential operator -- the fractional p-Laplacian. Nonlocal fractional operators,  in particular the fractional Laplacian,
have gained
a lot of popularity among
researchers working in a variety of fields. For instance, the fractional Laplacian  has been utilized
to model the dynamics in the Hamiltonian chaos in astrophysics (see \cite{10}), random motions such as the Brownian motion and the
Poisson process in physics (see \cite{11} and \cite{12}), the jump precess in finance and probability (see \cite{13}) as well as the
 the acoustic wave in mechanics. In the diffusion process, it has been used to derive heat kernel estimates for a large class of symmetric jump-type processes (see
 \cite{14}, \cite{15}).  The fractional Laplacian has  also been applied
 to the study of the game theory, image processing, L\'{e}vy processes, optimization and so on.  Readers who are interested in the application of the fractional Laplacians
 can refer to \cite{16}, \cite{17}, \cite{18} and the references therein.

The interest in the fractional operators continues to grow strong in this decade. More and more beautiful results,
 whose counterparts are powerful tools in elliptic PDE analysis,
have been proved in the fractional setting. The generalization, however,
 is no small feat due to the essential difference in how the fractional operators and the traditional differential operators are defined. In light of this, let's take a look at the
 fractional p-Laplacian. Let $s\in (0,1)$ and $p>1$.  The fractional p-Laplacian is defined as
\begin{eqnarray}\nonumber
(-\Delta)^s_p u(x)&=&C_{n,s,p}\lim_{ \varepsilon\rightarrow 0}\int_{R^n\backslash B_\varepsilon(x)}
\frac{|u(x)-u(y)|^{p-2}(u(x)-u(y))}{|x-y|^{n+ps}}dy\\
\label{20171041}
&=&C_{n,s,p}PV\int_{\R^n}\frac{|u(x)-u(y)|^{p-2}(u(x)-u(y))}{|x-y|^{n+ps}}dy,
\end{eqnarray}
where $PV$ stands for the Cauchy principal value.
To ensure that the integral in (\ref{20171041}) is well defined, we assume that
\[u\in C^{1,1}_{loc}(\Omega)\cap L_{sp}(R^n) \]
with \[L_{sp}=\{u\in L^1_{loc}(R^n) \mid \int_{\R^n}\frac{|u(x)|^{p-1}}{1+|x|^{n+sp}}dx<\infty\}.\]

When $s=1$ in Eq.(\ref{20171041}), it is the p-Lalacian.
When $p=2$,  Eq.(\ref{20171041}) becomes the nonlocal fractional Laplacian $(-\lap)^s$. A quick observation of the integral domain $R^n$ easily points to a characteristic shared among such integro-differential operators.
Different from the traditional differential operators, such as the Laplace operator, these are  not locally defined. To give  an example of what's new in nonlocal problems compared with local ones, we
consider the Dirichelet and the Neumann  problem on a bounded domain $\Omega\subset R^n$.
 To study the Laplacian problems, we require of information of solutions
  on the boundary $\partial\Omega$. But this is not enough for the fractional problems, which demand knowledge of solutions on both $\partial\Omega$ and $R^n\setminus\bar{\Omega}$. This
  raises a natural discussion in
  how to install appropriate boundary conditions in different cases
  so that the solutions can be extended in a way that preserve proper regularity in the whole space $R^n$. The challenge is especially true in computer-based simulations given that there is a limited amount of data we can gather over time. On top of this,
  when $p\neq2$,
the complexity increases because nonlinearity appears in the numerator.

In this paper, we are interested in the fractional p-Laplacian problems
with Dirichlet boundary conditions. Our first main result is a Hopf's lemma in a half-space. So far there are a few interesting results in the fractional settings  on the
Hopf's maximum principle.

In \cite{1}, Caffarelli et al. quoted a generalized Hopf's lemma for the smooth solution to a harmonic fractional equation
on a smooth domain $\Omega \subseteq R^n$. Either by
 the Harnack inequality or the Riesz potential, they claimed it true that
 if there is a
point $X_0\in \partial\Omega$  for which $v(X_0)=0$, then there exists
$\lambda> 0$ such that $v(x)\geq \lambda ((x-X_0)\cdot \nu (X_0))^\alpha$, where $\nu (X_0)$ is the inner normal to $\partial\Omega$ at $X_0$.

In \cite{8}, Greco and Servadei considered
 \begin{equation*}
(-\lap)^s u(x)\geq c(x)u(x), \quad x \in\Omega.
\end{equation*}
Assuming that $c(x)\leq 0$ in bounded domain $\Omega$, they derived that
\begin{equation}\label{201710251}
\inf \frac{u(x)}{(dist(x, \partial\Omega))^s}>0, \mbox{ as } x \ra \partial\Omega.
\end{equation}

Quite recently, Chen and Li\cite{5} proved a Hopf's lemma in terms of the boundary derivative
 for anti-symmetric functions on a half space through an elementary yet rather delicate analysis.
\begin{lem}[Chen-Li]
Assume that $w \in C^3_{loc}(\bar{\Sigma})$, $\overline{\underset{x \ra \partial\Sigma}{\lim}}=o(\frac{1}{[dist(x, \partial\Sigma)]^2})$, and
\begin{equation*}
\left\{\begin{array}{ll}
(-\lap)^s w(x)+c(x)w(x)=0, & in \;\Sigma,\\
w(x)>0,  & in \;\Sigma,\\
w(x^\la)=-w(x), &in\; \Sigma.
\end{array}
\right.
\end{equation*}
Then
\begin{equation*}
  \frac{\partial w}{\partial\nu}(x)<0, \quad x \in \partial\Sigma.
\end{equation*}
\end{lem}

In \cite{4}, Pezzo and Quass considered a fractional p-Laplacian problem
on a bounded domain $\Omega$  satisfying the interior ball condition:
\begin{equation}\label{201710181}
  (-\lap)^s_p u=c(x)|u|^{p-2}u, \quad x \in \Omega.
\end{equation}
 Under certain assumptions on $c(x)$, they obtained a similar result to that in
  (\ref{201710251}) for the weak super-solution of
 (\ref{201710181}).

Following the spirit in \cite{5}, we
 present a Hopf's lemma for $(-\lap)^s_p$ via the boundary derivative.

Let
\[T_\lambda=\{x\in \R^n \mid x_1=\lambda, \text{ for some } \lambda \in \R\}\]
be the moving planes,
\[\Sigma_\lambda=\{x\in \R^n \mid x_1>\lambda\}\]
be the region to the right of the plane $T_\lambda $,
\[ x^\lambda=(2\lambda-x_1,x_2,\cdots,x_n)\]
be the reflection of $x$ about  $T_\lambda$ and
\[w_\lambda(x)=u_\lambda(x)-u(x).\]

\begin{mthm}\label{201710250}
For $p\geq3$, assume that $u\in C^3_{loc}(\bar \Sigma)\cap L_{sp}$ and satisfies
\begin{equation}\label{2.2j}
\begin{aligned}
\begin{cases}
(-\Delta)_p^{s}u_\lambda(x)-(-\Delta)^{s}_p u(x)+c(x)w(x)=0, &\text{ in } \Sigma,\\
w(x)>0,&\text{ in } \Sigma,\\
w(x^\lambda)=-w(x),&\text{ in } \Sigma.
\end{cases}
\end{aligned}
\end{equation}
Let $\nu$ be the outward normal vector on $\partial\Sigma$. If \begin{equation}\label{2.1j}
\overline{\lim_{x\rightarrow \partial \Sigma}}c(x)=o(\frac{1}{[dist(x,\partial\Sigma)]^2}),
\end{equation}
then \begin{equation}\label{20171042}
\dfrac{\partial w}{\partial \nu}(x)<0, \quad x \in  \partial \Sigma.
\end{equation}
\end{mthm}

Following this we present our second main result-a Lipschitz boundary regularity for the fractional p-Laplacian.

In \cite{20}, Bogdan derived a boundary Harnack inequality for nonnegative
solutions for a harmonic fractional problem with Dirichlet condition.
Other boundary regularity for fractional equations were  obtained
by Caffarelli et al. in \cite{1} for a homogeneous fractional heat equation,
 and by Kim and Lee in \cite{22} for
a free boundary problem for the fractional Laplacian. In both papers the authors proved that the limit of $u(x)/dist(x, \partial\Omega)$ exists point-wise on the boundary.

 In a recent paper by Ros-Oton and Serra \cite{3}, the authors  considered
\begin{equation*}
\left\{\begin{array}{ll}
(-\lap)^s u=g,  &x \in \Omega,\\
u=0 , &x \in R^n\backslash\Omega.
\end{array}
\right.
\end{equation*}
Under the assumption that $g \in L^\infty (\Omega)$ for a bounded $\Omega$, they proved that the solution is $C^s(R^n)$ and
$\frac{u(x)}{dist(x, \partial\Omega)}$ is
$C^\alpha$ up to $\partial\Omega$ through a Krylov boundary Harnack inequality.

Later, in \cite{19} Chen et al. proved a similar result for the classical solutions through a good plain argument.
 The closest result to ours was obtained by  Iannizzotto et al. in \cite{23}. There the authors  proved $C^\alpha$ regularity, $\alpha \in (0,s]$, up to
the boundary for the weak solutions of a fractional p-Laplacian problem.
Their proof was carried out in the spirit of Krylov’s approach to boundary regularity and was quite complicated.

Inspired the work in \cite{19} and  \cite{23}, we
apply some of the ideas in \cite{19} on the following equation,
\begin{equation}\label{20171091}
\left\{\begin{array}{ll}
(-\Delta)^s_p u(x)=f(x),  &x \in \Omega,\\
u\equiv0, &x \in R^{n}\backslash\Omega.
\end{array}
\right.
\end{equation}

\begin{mthm}\label{20171092}
	Let $\Omega$ be a bounded domain in $R^n$ with exterior tangent spheres on the boundary, $s \in (0,1)$ and $p>2$. Assume that $\| f\|_{L^\infty(\Omega)}<\infty$ and $u(x) \in L_{ps}$.
If $u$ is a solutions of (\ref{20171091}), then there exists some $\nu \in (0, s)$ such that for $x$ close to the boundary
	\be
|u(x)|\leq c [dist(x, \partial\Sigma)]^\nu, \quad x \in \Sigma.
\ee
	\end{mthm}

\begin{mrem}
For $p=2$, $\nu$ can be up to $s$ (see \cite{3}).
\end{mrem}

For convenience's sake, we let
\[  |u(x)|^{p-2}u(x) =:[u(x)]^{p-1}.\]

Throughout the paper, we denote positive constants by $c$, $C_i$ whose values may vary from line to line.

\section{A Hopf's Lemma }

In this section, we prove Theorem \ref{201710250}.
For simplicity, in this section, we write $w_\lambda$ as $w$ and
$\Sigma_\la$ as $\Sigma$.

\textbf{Proof.}
To prove (\ref{20171042}), we develop a contradictive argument. Suppose there exists some $\tilde{x} \in \partial \Sigma$ such that (\ref{20171042}) is not true, then
\begin{equation}\label{20171043}
\dfrac{\partial w}{\partial \nu}(\tilde{x})=0.
\end{equation}

Without loss of generality, let $\lambda=0$ and $\tilde{x}$ be the origin. Let the ray from $\tilde{x}$ in the direction of $-\nu$ be the $x_1$ axis. By (\ref{20171043}) and the anti-symmetry of $w$, we know that $$\frac{\partial ^2 w}{\partial x_1^2}(0)=0.$$
For some $\bar{x}=(\bar{x}, 0') \in R^n$ and close to the origin, by the Taylor expansion, we obtain
\begin{equation}\label{2.6}
w(\bar{x})=w(0)+Dw(0)\cdot \bar{x}+\bar{x}\cdot D^2w(0)\cdot \bar{x}^T+O(|\bar{x}|^3)=O(|\bar{x}|^3).
\end{equation}
For simplicity's sake, let
 \begin{equation}
\delta= |\bar x_1|=dist(\bar{x},T_0).
\end{equation}
Then we have
$w(\bar{x})=O(\delta^3)$, and
\begin{equation}
|Dw(\bar{x})|=O(\delta^2),|D^2w(\bar{x})|=O(\delta).
\end{equation}
For $\bar{x}$ sufficiently close to the origin, i.e. $\delta$ sufficiently small, it's trivial that
\be\label{20171045}
c(\bar x)w(\bar{x})=o(1)\delta.
\ee
Using the estimate we have on $w_\lambda$ and its derivatives,  we can prove that for $\delta$ small and some $c_1>0$, it holds that
\be\label{20171044}
(-\Delta)_p^{s}u_\lambda(\bar x)-(-\Delta)^{s}_p u(\bar x)\leq - \frac{c_1}{4}\delta.
\ee
We postpone the proof of (\ref{20171044}) for the moment. Combining
(\ref{20171045}) and
(\ref{20171044}), we arrive at
$$(-\Delta)_p^{s}u_\lambda(\bar x)-(-\Delta)^{s}_p u(\bar x)+c(\bar x)w(\bar{x})<0.$$
This contradicts to (\ref{2.2j}) and thus proves the theorem.

Now we prove (\ref{20171044}).

Recall that  $y^\la=y^0=(-y_1,y')$. By (\ref{20171041}), we have
\begin{eqnarray}\nonumber
&& (-\Delta)^s_p u_\lambda(\bar x) -(-\Delta)^s_p u(\bar x)\\\nonumber
&= &C_{n,s,p}PV\int_{R^n}\dfrac{(u_\lambda(\bar x)-u_\lambda(y))^{p-1}-(u(\bar x)-u(y))^{p-1}}{|\bar x-y|^{n+ps}}dy\\\nonumber
&= &C_{n,s,p}PV\int_{\Sigma}\dfrac{(u_\lambda(\bar x)-u_\lambda(y))^{p-1}-(u(\bar x)-u(y))^{p-1}}{|\bar x-y|^{n+ps}}dy\\\nonumber
&&\qquad
+C_{n,s,p}\int_{R^n/\Sigma}\dfrac{(u_\lambda(\bar x)-u_\lambda(y))^{p-1}-(u(\bar x)-u(y))^{p-1}}{|\bar x-y|^{n+ps}}dy\\\nonumber
&=&C_{n,s,p}PV\int_{\Sigma}\dfrac{(u_\lambda(\bar x)-u_\lambda(y))^{p-1}-(u(\bar x)-u(y))^{p-1}}{|\bar x-y|^{n+ps}}dy\\\nonumber
&&\qquad
+C_{n,s,p}\int_{\Sigma}\dfrac{(u_\lambda(\bar x)-u(y))^{p-1}-(u(\bar x)-u_\lambda(y))^{p-1}}{|\bar x-y^0|^{n+ps}}dy\\\nonumber
&=&C_{n,s,p}PV\int_{\Sigma}\Bigl((u_\lambda(\bar x)-u_\lambda(y))^{p-1}-(u(\bar x)-u(y))^{p-1}\Bigl)\\\nonumber
&&\cdot \Bigl (\dfrac{1}{|\bar x-y|^{n+ps}}-\dfrac{1}{|\bar x-y^0|^{n+ps}}\Bigl )dy\\\nonumber
&&+C_{n,s,p}\int_{\Sigma}\dfrac{(u_\lambda(\bar x)-u_\lambda(y))^{p-1}-(u(\bar x)-u(y))^{p-1}}{|\bar x-y^0|^{n+ps}}\\\nonumber
&&\frac{+(u_\lambda(\bar x)-u(y))^{p-1}-(u(\bar x)-u_\lambda(y))^{p-1}}{}dy
\\\label{2.1}
&=:&C_{n,s,p}PV\int_{\Sigma}I\,dy+C_{n,s,p}\int_{\Sigma}II\,dy.
\end{eqnarray}

We first take care of $\int_{\Sigma}IIdy$ for later.

Let $R_o>0$ be a given positive number. Then
$$\Sigma= \big(\Sigma\cap B_{R_o}(\bar{x}) \big) \cup
 \big( \Sigma\cap B^c_{R_o}(\bar{x}) \big).$$

 For
 $y\in \Sigma \cap B_{R_o}(\bar{x})$, by the mean value theorem we have
\begin{eqnarray}\nonumber
&& (u_\lambda(\bar x)-u_\lambda(y))^{p-1}-(u(\bar x)-u(y))^{p-1}+(u_\lambda(\bar x)-u(y))^{p-1}-(u(\bar x)-u_\lambda(y))^{p-1}\\\nonumber
&=&(u_\lambda(\bar x)-u_\lambda(y))^{p-1}-(u(\bar x)-u_\lambda(y))^{p-1}+(u_\lambda(\bar x)-u(y))^{p-1}-(u(\bar x)-u(y))^{p-1}\\\nonumber
&=&(p-1)(|\xi_1|^{p-2} +|\xi_2|^{p-2}  )w_\lambda(\bar x) \\\label{2.2}
&\leq &c w_\lambda(\bar x) |\bar x-y^0| ^{p-2},                                                                                                                                                                                                                                                                                                                                                                                                                                        \end{eqnarray}
with $\xi_1$ between $u_\lambda(\bar x)-u_\lambda(y)$ and $u(\bar x)-u_\lambda(y)$, $\xi_2$   between $u_\lambda(\bar x)-u(y)$ and $u(\bar x)-u(y)$.  The last inequality holds because, under the assumption $w(y)>0$ for
$y \in \Sigma$, we have
\begin{eqnarray*}
|\xi_1|&\leq & \max\{|u_\lambda(\bar x)-u_\lambda(y)|, |u(\bar x)-u_\lambda(y)|\}\\
&<&\max\{|u_\lambda(\bar x)-u_\lambda(y)|, |u(\bar x)-u(y)|\}\\
&\leq &c\max\{|\bar x-y^0|, |\bar x-y|\}=c|\bar x-y^0|.
\end{eqnarray*}
Similarly, one can show that
$$|\xi_2|< c\max\{|\bar x-y^0|, |\bar x-y|\}=c|\bar x-y^0|.$$

From (\ref{2.2}), for $\delta$ sufficiently small we deduce that
\begin{eqnarray}\nonumber
 \int_{\Sigma \cap B_{R_o}(\bar{x})}|II|dy&=& \int_{\Sigma \cap B_{R_o}(\bar{x})} IIdy\\\nonumber
&\leq& c\int_{\Sigma \cap B_{R_o}(\bar{x})}\dfrac{w(\bar x)}{|\bar x-y^0|^{n+ps-p+2}}dy\\\nonumber
&\leq& c w(\bar{x})\int_{ B_{2R_o}(\bar{x})\backslash B_\delta(\bar{x}) }\dfrac{1}{|\bar x-y^0|^{n+ps-p+2}}dy\\\label{2.10j}
 &\leq&  c_1 \max\{\delta^{1+p-ps},\delta^2\}.
\end{eqnarray}

For $y\in \Sigma \cap B^c_{R_o}(\bar{x})$, using $u \in L_{sp}$ and the
H\"{o}lder inequality we have
\begin{eqnarray*}
&&\int_{\Sigma \cap B^c_{R_o}(\bar{x})} II \, dy\\
 & \leq& cw(\bar x)\int_{\Sigma \cap B^c_{R_o}(\bar{x}) } \dfrac{|u(\bar x)|^{p-2}+|u_\lambda(\bar x)|^{p-2}+|u(y)|^{p-2}+|u_\lambda(y)|^{p-2}}{|\bar x-y^0|^{n+ps}}dy\\
&\leq & c\,w(\bar x)\Big[C\int_{ |y|\geq R/2}\dfrac{1}{(1+|y|)^{n+ps}}dy+2\int_{ |y|\geq R/2}\dfrac{|u(y)|^{p-2}}{(1+|y|)^{n+ps}}dy\Big]\\
&\leq& c\,w(\bar x)\bigg(C+2\big(\int_{ |y|\geq R/2}\dfrac{|u(y)|^{p-1}}{(1+|y|)^{n+ps}}dy\big)^{\frac{p-2}{p-1}}
\big(\int_{ |y|\geq R/2}\dfrac{1}{(1+|y|)^{n+ps}}dy\big)^{\frac{1}{p-1}}\bigg)\\
&\leq& c\,w(\bar x)\\
&\leq& c\,\delta^3.
\end{eqnarray*}
Together with (\ref{2.10j}), it shows that for $\delta$ small we have
\begin{equation}\label{2.13j}
\int_{\Sigma}IIdy\leq  c \max\{\delta^{1+p-ps},\delta^2\}.
\end{equation}

Next we estimate $\int_{\Sigma}Idy$.

\begin{figure}
\centering
\includegraphics{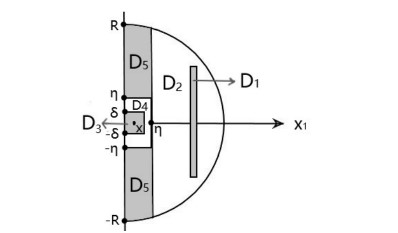}
\label{p1}
\caption{Subregions}
\end{figure}

For some $R>>1$ large, let  $B^+_R(0)= \{x \in B_R(0) \mid x_1>0  \}$.
To take care of the possible singularities, we divide $B^+_R(0)$ into five subregions( see Fig.\ref{p1}) defined as below.
\begin{equation*}
D_1=\{x\mid 1\leq x_1\leq 2,\,|x'|\leq 1\},
\end{equation*}
\begin{equation*}
D_2=\{x \in B_R(0)\mid  x_1\geq \eta\},
\end{equation*}
\begin{equation*}
D_3=\{x \mid 0\leq x_1\leq 2\delta, \, |x'|<\delta\},
\end{equation*}
\begin{equation*}
 D_4=\{x \mid 0\leq x_1\leq \eta, \: |x'|<\eta, \:x \not\in D_3\},
\end{equation*}
\begin{equation*}
 D_5=\{x  \in B_R(0)\mid 0 \leq x_1\leq \eta,\: |x'|>\eta\}.
\end{equation*}

We estimate the integral in each region accordingly.
 Later, we will discuss the requirements that $R$ and $\eta$ must satisfy.
 Roughly speaking, we need to take $R$ sufficiently large and $\eta> \delta$ sufficiently small.

We start with $D_1$.
By the mean value theorem we have
\begin{eqnarray}\nonumber
&&\dfrac{1}{|\bar x-y|^{n+ps}}-\dfrac{1}{|\bar x-y^0|^{n+ps}}\\\nonumber
&=&(-\frac{n+ps}{2})\frac{1}{|\xi_3|^{\frac{n+ps}{2}+1}}\bigl(|\bar x-y|^2-|\bar x-y^0|^2\bigl )\\\label{2.11}
&=&\frac{n+ps}{2} \frac{1}{|\xi_3|^{\frac{n+ps}{2}+1}}4\bar x_1y_1
\end{eqnarray}
with
$$|\bar x-y|^2\leq \xi_3\leq |\bar x-y^0|^2,$$
and
\be\label{2.12}
(u_\lambda(\bar x)-u_\lambda(y))^{p-1}-(u(\bar x)-u(y))^{p-1}\\
=(p-1)|\xi_4|^{p-2}[w(\bar x)-w(y)],
\ee
where $\xi_4$ is between $u_\lambda(\bar x)-u_\lambda(y)$ and $u(\bar x)-u(y)$.

Since $w(x)>0$ in $\Sigma$ and $w(0)=0$, for  $y\in D_1$ and
 $\bar x$ sufficiently close to the origin, it is trivial that
\be\label{2.24jjj}
w(\bar x)-w_\lambda(y)<-c<0.
\ee
Hence
$$u_\lambda(\bar x)-u_\lambda(y)<u(\bar x)-u(y).$$
Together with (\ref{2.12}), it shows that
$$\xi_4 \neq 0, \quad y \in \Sigma.$$
Therefore there exists some $c$ such that
$$|\xi_4|\geq c>0.$$
Combine this result with (\ref{2.11}) and (\ref{2.12}), it gives
\begin{eqnarray}\nonumber
  &&\int_{D_1}\Bigl([u_\lambda(\bar x)-u_\lambda(y)]^{p-1}-[u(\bar x)-u(y)]^{p-1}\Bigl)\Bigl (\dfrac{1}{|\bar x-y|^{n+ps}}-\dfrac{1}{|\bar x-y^0|^{n+ps}}\Bigl )dy\\\nonumber
&\leq& c\int_{D_1}|\xi_4|^{p-2}[w(\bar x)-w(y)]\frac{x_1 y_1}{|\bar x-y^0|^{n+ps+2}}dy\\\label{20171081}
&\leq&  -\int_{D_1} c\delta dx\leq -c_1\delta.
\end{eqnarray}

We estimate the integral on $D_2$. Later, in the proof for
 $D_4$ and $D_5$, we will discuss the ranges of $\eta$ and $R$ respectively. For now, we assume both $R$ and $\eta$ have already been selected and fixed. Then it's obvious that $$w_\lambda(x)-w_\lambda(y)\leq 0, \quad y\in \Omega_{R,\eta}, as \; \delta \ra 0.$$ 
Thus
\begin{eqnarray}\nonumber
   && \int_{D_2 }\Bigl([u_\lambda(\bar x)-u_\lambda(y)]^{p-1}-[u(\bar x)-u(y)]^{p-1}\Bigl)\Bigl (\dfrac{1}{|\bar x-y|^{n+ps}}-\dfrac{1}{|\bar x-y^0|^{n+ps}}\Bigl )dy \\\nonumber
 & \leq &\int_{D_1}\Bigl(u_\lambda(\bar x)-u_\lambda(y))^{p-1}-(u(\bar x)-u(y))^{p-1}\Bigl)\Bigl (\dfrac{1}{|\bar x-y|^{n+ps}}-\dfrac{1}{|\bar x-y^0|^{n+ps}}\Bigl )dy \\\label{2.27}
   &\leq &-c_1\delta
\end{eqnarray}

On $D_3$, we separate the integrand $I$ into two pieces. On one hand,
by Taylor expansion, we have
\begin{eqnarray}\nonumber
  &&\bigg|\int_{D_3}\frac{[u_\lambda(\bar x)-u_\lambda(y)]^{p-1}-[u(\bar x)-u(y)]^{p-1}}{|\bar x-y^0|^{n+ps}}dy\bigg|\\\nonumber
&\leq& c\int_{D_3}\frac{|\xi_4|^{p-2}}{|\bar x-y^0|^{n+ps}}\big|w(\bar x)-w(y)\big|dy
\\\nonumber
&\leq&c\int_{D_3}\frac{|\bar x-y|^{p-2}}{|\bar x-y^0|^{n+ps}}
\big(|Dw(\bar x)\cdot(\bar x-y)|+ |(\bar x-y)\cdot D^2 w(\bar x)\cdot(\bar x-y)^{T}|\\\nonumber
&&\qquad+|O(|\bar x-y|^3)|\big)dy\\\nonumber
&\leq&c\delta^2\int_{D_3}\frac{1}{|\bar x-y^0|^{n+ps-p+1}}dy\\\label{20171072}
&\leq&c \max\{\delta^2, \delta^{1+p-ps}\}.
\end{eqnarray}
The last inequality is true as $\delta \ra 0$.

On the other hand, for $ \xi_5$ between $u_\la(\bar x)-u_\lambda(y)$ and $Du_{\la}(\bar x)\cdot(\bar x-y)$, we have
\begin{eqnarray*}
   && \int_{D_3}\frac{[u_\lambda(\bar x)-u_\lambda(y)]^{p-1}}
   {|\bar x-y|^{n+ps}}dy\\
   &=&\int_{D_3}\frac{[Du_{\la}(\bar x)\cdot(\bar x-y)]^{p-1}+ (p-1)|\xi_5|^{p-2}[(\bar x-y)\cdot D^2 u_{\la}(\bar x)\cdot(\bar x-y)^{T}}{|\bar x-y|^{n+ps}} \\
   &&\frac{+O(|\bar x-y|^3)]}{}dy\\
   &=& \int_{D_3}\frac{(p-1)|\xi_5|^{p-2}[(\bar x-y)\cdot D^2 u_{\la}(\bar x)\cdot(\bar x-y)^{T}+O(|\bar x-y|^3)]}{|\bar x-y|^{n+ps}}dy\\
   &\leq & \int_{D_3}(p-1)\frac{|\xi_5|^{p-2}(\bar x-y)\cdot D^2 u_{\la}(\bar x)\cdot(\bar x-y)^{T}}{|\bar x-y|^{n+ps}}dy+ O(1)\delta^{1+p-ps}.
\end{eqnarray*}
We obtain the second to last equation from the fact that
$$ \int_{D_3}\frac{[Du_{\la}(\bar x)\cdot(\bar x-y)]^{p-1}}{|\bar x-y|^{n+ps}}dy=0,$$
as a result of the symmetry of $D_3$ with respect to $\bar x$. Similarly,
for $ \xi_6$ between $u(\bar x)-u(y)$ and $Du(\bar x)\cdot(\bar x-y)$,
 we have
\begin{eqnarray*}
   && \int_{D_3}\frac{[u(\bar x)-u(y)]^{p-1}}
   {|\bar x-y|^{n+ps}}dy\\
&\leq&\int_{D_3}(p-1)\frac{|\xi_6|^{p-2}(\bar x-y)\cdot D^2 u(\bar x)\cdot(\bar x-y)^{T}}{|\bar x-y|^{n+ps}}dy+O(1)\delta^{1+p-ps}.
\end{eqnarray*}
Therefore, it follows that
\begin{eqnarray}\nonumber
  &&\bigg|\int_{D_3} \dfrac{[u_\lambda(\bar x)-u_\lambda(y)]^{p-1}-[u(\bar x)-u(y)]^{p-1}}
  {|\bar x-y|^{n+ps}}dy\bigg|\\\nonumber
  &\leq & c(p-1)\bigg|\int_{D_3} \dfrac{(\bar x-y)\cdot[ |\xi_5|^{p-2} D^2 u_{\la}(\bar x)-|\xi_6|^{p-2} D^2 u(\bar x)]
  \cdot(\bar x-y)^{T}}
  {|\bar x-y|^{n+ps}}dy\bigg|\\\nonumber
  &&+c\delta^{1+p-ps}\\\nonumber
  &=&\bigg|\int_{D_3} \dfrac{(\bar x-y)\cdot |\xi_5|^{p-2} D^2 w(\bar x)\cdot(\bar x-y)^{T}}
  {|\bar x-y|^{n+ps}}\\\nonumber
  &&\frac{+(|\xi_5|^{p-2}-|\xi_6|^{p-2})(\bar x-y)\cdot D^2 u(\bar x)\cdot(\bar x-y)^{T}}{}dy\bigg|+c\delta^{1+p-ps}\\\label{20171074}
&\leq& c\delta^{1+p-ps}.
\end{eqnarray}
Combining (\ref{20171072}) with (\ref{20171074}) it gives
\be\label{20171071}
|\int_{D_3} I\,dy| \leq c \max\{\delta^2, \delta^{1+p-ps}\}.
\ee

Below we deal with $D_4$.
\begin{eqnarray}\nonumber
  &&\bigg|\int_{D_4}\Bigl([u_\lambda(\bar x)-u_\lambda(y)]^{p-1}-[u(\bar x)-u(y)]^{p-1}\Bigl)\Bigl (\dfrac{1}{|\bar x-y|^{n+ps}}-\dfrac{1}{|\bar x-y^0|^{n+ps}}\Bigl )dy\bigg|\\\nonumber
&\leq& c\int_{D_4}|\xi_4|^{p-2}|w(\bar x)-w(y)\big|\frac{x_1 y_1}{|\bar x-y^0|^{n+ps+2}}dy\\\nonumber
&\leq& c\int_{D_4}|\bar x-y|^{p-2}
\big(|Dw(\bar x)\cdot(\bar x-y)|+ |(\bar x-y)\cdot D^2 w(\bar x)\cdot(\bar x-y)^{T}|\\\nonumber
&&\qquad+|O(|\bar x-y|^3)|\big)\frac{\delta}{|\bar x-y|^{n+ps+1}}dy\\\nonumber
&\leq&  c\delta\int_{B_{2\eta}(\bar x)\backslash B_{\delta}(\bar x )}
\frac{1}{|\bar x-y|^{n+ps-p}} dy\\\nonumber
&=&  c\delta \frac{(2\eta)^{p-ps}-\delta^{p-ps}}{p-ps}\\\label{20171082}
&\leq& \frac{ c_1}{8}\delta.
\end{eqnarray}
The last inequality is true when $\eta$ is sufficiently small.

On $D_5$, we have
\begin{eqnarray}\nonumber
  &&\bigg|\int_{D_5}\Bigl([u_\lambda(\bar x)-u_\lambda(y)]^{p-1}-[u(\bar x)-u(y)]^{p-1}\Bigl)\Bigl (\dfrac{1}{|\bar x-y|^{n+ps}}-\dfrac{1}{|\bar x-y^0|^{n+ps}}\Bigl )dy\bigg|\\\nonumber
&\leq& c\int_{D_5}|\xi_4|^{p-2}|w(\bar x)-w(y)\big|\frac{x_1 y_1}{|\bar x-y^0|^{n+ps+2}}dy\\\nonumber
&\leq& c\int_{D_5}|\bar x-y|^{p-2}
\big(|Dw(\bar x)\cdot(\bar x-y)|+ |(\bar x-y)\cdot D^2 w(\bar x)\cdot(\bar x-y)^{T}|\\\nonumber
&&\qquad+|O(|\bar x-y|^3)|\big)\frac{\delta \eta}{|\bar x-y|^{n+ps+2}}dy\\\nonumber
&\leq&  c\delta \eta\int_{B_{2R}(\bar x)\backslash B_{\eta}(\bar x)}
\frac{1}{|\bar x-y|^{n+ps-p+1}} dy\\\nonumber
&=&  c\delta \eta \frac{(2R)^{3+p-ps}-(\eta)^{3+p-ps}}{3+p-ps}\\\label{20171083}
&\leq& \frac{ c_1}{8}\delta.
\end{eqnarray}
The validity of the last inequality results from $\eta$
being sufficiently small for $R$ fixed.

Gathering the estimates on $D_i$, $i=1,2,3,4,5$, that is, (\ref{20171081}), (\ref{2.27}),  (\ref{20171071}),  (\ref{20171082}) and  (\ref{20171083}), it
shows that for $\eta$ sufficiently small,
\be\label{20171087}
\int_{B_R^+(0)}I \,dy\leq -c_1\delta.
\ee

What remains to do is the integral on $\Sigma\backslash B^+_R(0) $.
\begin{eqnarray}\nonumber
 &&\bigl|\int_{\Sigma\backslash B^+_R(0)}
 \Bigl([u_\lambda(\bar x)-u_\lambda(y)]^{p-1}-[u(\bar x)-u(y)]^{p-1}\Bigl)
 \Bigl (\dfrac{1}{|\bar x-y|^{n+ps}}-\dfrac{1}{|\bar x-y^0|^{n+ps}}\Bigl )dy\bigl |\\\nonumber
& \leq &c\delta\int_{\Sigma\backslash B^+_R(0)} \frac{|u(\bar x)|^{p-1}+|u_\lambda(\bar x)|^{p-1}+|u_\lambda (y)|^{p-1}+|u(y)|^{p-1}}{|\bar x-y|^{n+ps+1}}dy\\ \nonumber
&\leq & c\delta\int_{\Sigma\backslash B^+_R(0)} \frac{|u(\bar x)|^{p-1}+|u_\lambda(\bar x)|^{p-1}}{|\bar x-y|^{n+ps+1}}dy
+ \frac{c\delta}{R}\int_{\R^n} \frac{ |u(y)|^{p-1}}{(1+|y|)^{n+ps}}dy\\\nonumber
&\leq & \frac{c\delta}{R^{1+ps}}+\frac{c\delta}{R}\\\label{20171086}
&\leq &\frac{c_1}{8}\delta.
\end{eqnarray}
The last inequality holds when $R$ is sufficiently large. Together with (\ref{20171087}), it gives
\be\label{20171088}
\int_{\Sigma}I \,dy\leq -\frac{c_1\delta}{2}.
\ee

Combining this with (\ref{2.13j}), for $\delta$ sufficiently small, we conclude that
$$
(-\Delta)^s_p u_\lambda(\bar x) -(-\Delta)^s_p u(\bar x)\leq -\frac{c_1\delta}{4}.
$$
This proves (\ref{20171044}) and completes the proof of the theorem.


\section{Boundary Regularity}

In this section we prove Theorem \ref{20171092}.
Here the analysis of regularity up to the boundary is based on the existence of some super-solution, sometimes referred to as the barrier function
  in boundary regularity analysis, to the fractional p-Laplacian equation.
To construct the barrier function, we begin with an equation in $R^+_1:=\{x \in R \mid x>0   \}$, whose solution is known explicitly.
 \begin{lem}\label{l3.1}
 	For $0<\nu<s$,
 	\begin{equation}
 	(-\Delta)^s_p(x^\nu_+)=C_\nu x_+^{(p-1)\nu-ps},\quad x\in \R^+,
 	\end{equation}
 with \[C_\nu=\int^{+\infty}_{-\infty}\frac{(1-z^\nu_+)^{p-1}}{|1-z|^{1+ps}}dz>0\].
 \end{lem}

\textbf{Proof.}
 Since $x>0$, we have $x_+=x$, and
 \begin{equation}
 \begin{aligned}
 (-\Delta)^s_p(x^\nu_+)&=\int_{\R}\frac{(x^\nu_+-y^\nu_+)^{p-1}}{|x-y|^{1+ps}}dy\\
 &=\int_{-\infty}^{+\infty}\frac{x^{(p-1)\nu}(1-z^\nu_+)^{p-1}}{x^{1+ps}|1-z|^{1+ps}}xdz \;( y=xz) \\
 &=x^{(p-1)\nu-ps}\int_{-\infty}^{+\infty}\frac{(1-z^\nu_+)^{p-1}}{|1-z|^{1+ps}}dz\\
 &=C_\nu x_+^{(p-1)\nu-ps},
 \end{aligned}
 \end{equation}
with $\displaystyle C_\nu=\int^{+\infty}_{-\infty}\dfrac{(1-z^\nu_+)^{p-1}}{|1-z|^{1+ps}}dz$.
Then
 \begin{eqnarray}\nonumber
 C_\nu&=&\int^{+\infty}_{0}\dfrac{(1-z^\nu_+)^{p-1}}{|1-z|^{1+ps}}dz+\int^{0}_{-\infty}\dfrac{(1-z^\nu_+)^{p-1}}{|1-z|^{1+ps}}dz\\\label{3.4}
 &=& \int^{+\infty}_{0}\dfrac{(1-z^\nu_+)^{p-1}}{|1-z|^{1+ps}}dz+\frac{1}{ps}.
 \end{eqnarray}

For $0<\nu\leq\frac{ps-1}{p-1}$,
 \begin{equation*}
 \begin{aligned}
 &\ \ \ \int^{+\infty}_{0}\dfrac{(1-z^\nu_+)^{p-1}}{|1-z|^{1+ps}}dz
 \\&=\int^{1}_{0}\dfrac{(1-z^\nu)^{p-1}}{|1-z|^{1+ps}}dz+\int^{+\infty}_{1}\dfrac{(1-z^\nu)^{p-1}}{|1-z|^{1+ps}}dz\\
 &=\int^{1}_{0}\dfrac{(1-z^\nu)^{p-1}}{|1-z|^{1+ps}}dz+\int^{1}_{0}\dfrac{(1-w^{-\nu})^{p-1}}{|1-w^{-1}|^{1+ps}}\frac{1}{w^2}dw\\
 &=\int^{1}_{0}\dfrac{(1-z^\nu)^{p-1}}{|1-z|^{1+ps}}dz+\int^{1}_{0}-\dfrac{(1-w^\nu)^{p-1}}{|1-w|^{1+ps}}w^{1+ps-\nu(p-1)-2}dz\\
 &=\int^{1}_{0}\dfrac{(1-z^\nu)^{p-1}}{|1-z|^{1+ps}}(1-z^{ps-\nu(p-1)-1})dz\geq0.
 \end{aligned}
 \end{equation*}
Together with (\ref{3.4}) it implies that
\begin{equation}\label{3.6}
C_\nu>0,  \quad  \mbox{for }0<\nu<\frac{ps-1}{p-1}.
\end{equation}

 To continue, we need Lemma 3.1 in \cite{23} which states that
 \begin{equation}
 (-\Delta)^s_p(x_+^s)|_{x=1}=0, \quad x\in R.
 \end{equation}
 Then  for  $\frac{ps-1}{p-1}<\nu<s$, it follows that
 \begin{eqnarray*}\nonumber
   C_\nu &=&  \int^{+\infty}_{-\infty}\dfrac{(1-z^\nu_+)^{p-1}-(1-z^s_+)^{p-1}}{|1-z|^{1+ps}}dz\\\nonumber
 &=&\int^{+\infty}_{0}\dfrac{(1-z^\nu_+)^{p-1}-(1-z^s_+)^{p-1}}{|1-z|^{1+ps}}dz\\ \nonumber &=&\int^{+\infty}_1\dfrac{(1-z^\nu_+)^{p-1}-(1-z^s_+)^{p-1}}{|1-z|^{1+ps}}dz
 +\int^1_0\dfrac{(1-z^\nu_+)^{p-1}-(1-z^s_+)^{p-1}}{|1-z|^{1+ps}}dz\\\nonumber
&=&\int^{1}_0\dfrac{(1-z^{-\nu}_+)^{p-1}-(1-z^{-s}_+)^{p-1}}{|1-z^{-1}|^{1+ps}}\frac{1}{z^2}dz
+\int^1_0\dfrac{(1-z^\nu_+)^{p-1}-(1-z^s_+)^{p-1}}{|1-z|^{1+ps}}dz\\\nonumber
&=&\int^{1}_0\dfrac{-(1-z^{\nu})^{p-1}z^{ps-1-(p-1)\nu}+(1-z^{s}_+)^{p-1}z^{s-1}}{|1-z|^{1+ps}}\frac{1}{z^2}dz
\\&& +\int^1_0\dfrac{(1-z^\nu_+)^{p-1}-(1-z^s)^{p-1}}{|1-z|^{1+ps}}dz\\\nonumber
&=&\int^1_0\dfrac{(1-z^\nu)^{p-1}-(1-z^s)^{p-1}}{|1-z|^{1+ps}}(1-z^{ps-1-(p-1)\nu})dz\\
&& +\int^1_0\dfrac{(1-z^s)^{p-1}}{|1-z|^{1+ps}}(z^{s-1}-z^{ps-1-(p-1)\nu})dz>0.
 \end{eqnarray*}
Together with (\ref{3.6}), we conclude that
$$
   C_{\nu}>0 \text{ for } \nu\in (0, s).
$$

Next we generalize Lemma \ref{l3.1} to $n-$dimensions.
Let $R^+_n:=\{x \in R^n \mid x_n>0   \}$.

\begin{cor}\label{c3.1}
For $0<\nu<s$,
		\begin{equation}
		(-\Delta)^s_p(x_n)^\nu_+=C_{\nu,n} (x_n)_+^{(p-1)\nu-ps},\text{ for } x \in R^n_+,
		\end{equation}
		with \[C_{\nu,n}=C_{\nu}\int_0^{\infty}\frac{t^{n-2}}{(1+t^2)^{\frac{n+ps}{2}}}dt >0.\]
\end{cor}

\textbf{Proof.}
 Let $x=(x',x_n)\in \R^n,r=|x'-y'| \text{ and }\tau=|x_n-y_n|$.
From Lemma \ref{l3.1} we have
\begin{eqnarray*}
(-\Delta)^s_p({x_n})^\nu_+&=&\int_{\R^n}\frac{((x_n)^\nu_+-{(y_n)^\nu_+})^{p-1}}{|x-y|^{n+ps}}dy\\
&=&\int_{\R^n}\frac{({(x_n)^\nu_+}-{(y_n)^\nu_+})^{p-1}}{|(x_n-y_n)^2+(x'-y')^2|^{\frac{n+ps}{2}}}dy\\
&=&\int_{-\infty}^{+\infty}({(x_n)^\nu_+}-{(y_n)^\nu_+})^{p-1}
\int_0^{\infty}\frac{w_{n-2}r^{n-2}}{|\tau^2+r^2|^{\frac{n+ps}{2}}}dr\,dy_n\\
&=&\int_{-\infty}^{+\infty}
\frac{({(x_n)^\nu_+}-{(y_n)^\nu_+})^{p-1}}{|x_n-y_n|^{1+ps}}dy_n
\int_0^{\infty}\frac{w_{n-2}t^{n-2}}{|1+t^2|^{\frac{n+ps}{2}}}dt\: (r=\tau t)\\
&=&(x_n)^{\nu(p-1)-ps}_+ C_\nu \int_0^{\infty}\frac{w_{n-2}t^{n-2}}{|1+t^2|^{\frac{n+ps}{2}}}dt\\
&:=&(x_n)^{\nu(p-1)-ps}_+ C_{\nu,n}.
\end{eqnarray*}


Now we are ready to construct the barrier function.
\begin{lem}\label{20171098}
	Let $\phi(x)=(|x|^2-1)^\nu_+$ in $R^n$ with $\nu\in (0,s)$. Then there exists some  $\epsilon>0$ small and  $C_0>0$ such that
	\begin{equation}\label{20171093}
	(-\Delta)^s_p\phi(x)\geq C_0(|x|-1)^{\nu(p-1)-ps}, \quad x \in B_{1+\epsilon}(0)\backslash B_1(0).
	\end{equation}
\end{lem}

\textbf{Proof.}
 To prove the lemma, we argue by contradiction. Suppose (\ref{20171093}) is not true, then there
 exists a sequence $\{x^k\} \in B_1(0)$ so that
 $|x^k|\rightarrow 1$ and
 \be\label{20171094}
 (-\Delta)^s_p\phi(x^k)(|x^k|-1)^{ps-\nu(p-1)} \ra 0,  \mbox{ as } \quad  k \ra \infty.
 \ee

Without loss of generality, let  $x^k=(0,1+d_k)$. Then
$$d_k=|x^k|-1 \ra 0, \quad \mbox{as } k \ra \infty.$$
Here we use an equivalent form of (\ref{20171041}) via the difference quotient
$$(-\Delta)_p^s \phi(x^k)=\frac{C_{n,s,p}}{2}\int_{\R^n}\frac{[\phi(x^k)-\phi(x^k+y)]^{p-1}+[\phi(x^k)-\phi(x^k-y)]^{p-1}}{|y|^{n+ps}}dy.$$
Then by Lemma \ref{c3.1}, we have
 \begin{equation}
 \begin{aligned}
 &\ \ \ \ (|x^k-1|)^{ps-\nu(p-1)}(-\Delta)_p^s \phi(x^k)\\
 &=\frac{d_k^{ps-\nu(p-1)}}{2}C_{n,s,p}\int_{\R^n}
 \frac{[\phi(x^k)-\phi(x^k+y)]^{p-1}+[\phi(x^k)-\phi(x^k-y)]^{p-1}}{|y|^{n+ps}}dy\\
 &=C_{n,s,p}d_k^{ps-\nu(p-1)}\bigl(\int_{\R^n}\frac{\Bigl[(|x^k|^2-1)^\nu_+-(|x^k+y|^2-1)^\nu_+\Bigl]^{p-1}}{|y|^{n+ps}}dy\\
 &\ \ +\int_{\R^n}\frac{\Bigl[(|x^k|^2-1)^\nu_+-(|x^k-y|^2-1)^\nu_+\Bigl]^{p-1}}{|y|^{n+ps}}dy\bigl)\\
 &=C_{n,s,p}d_k^{ps-\nu(p-1)}\bigl(\int_{\R^n}\frac{\big[(d_k^2+2d_k)^\nu_+-(d_k^2+2d_k +2(1+d_k)y_n+|y|^2)^\nu_+\big]^{p-1}}{|y|^{n+ps}}\\
 &\ \ +\frac{\big[(d_k^2+2d_k)^\nu_+-(d_k^2+2d_k -2(1+d_k)y_n+|y|^2)_+^\nu\big]^{p-1}}{}dy\bigl)\\
 &=\frac{C_{n,s,p}}{2}\int_{\R^n}\bigg(\frac{[(d_k+2)^\nu_+-(d_k+2 +2(1+d_k)z_n+d_k|z|^2)^\nu_+]^{p-1}}{|z|^{n+ps}}\\
 &\ \ \frac{+\big[(d_k+2)^\nu_+-(d_k+2 -2(1+d_k)z_n+d_k|z|^2)_+^\nu\big]^{p-1}}{}\bigg)dz \quad  (y=d_k z )\\
 &\ \ \ \rightarrow\frac{C_{n,s,p}}{2}
 \int_{\R^n}\frac{(2^\nu-(2+2z_n)_+^\nu)^{p-1}+(2^\nu-(2-2z_n)_+^\nu)^{p-1}}{|z|^{n+ps}}dz\\
 &=2^{(p-1)\nu-1}C_{n,s,p}\int_{\R^n}\frac{(1-(1+z_n)_+^\nu)^{p-1}+(1-(1-z_n)^\nu_+)^{p-1}}{|z|^{n+ps}}dz\\
 &=2^{(p-1)\nu}(-\Delta)_p^s(x_n)^\nu_+|_{x_n=1}\\
 &=2^{(p-1)\nu}C_{\nu,n}>0.
 \end{aligned}
 \end{equation}
This is a contradiction with (\ref{20171094}).

In addition to the barrier function, we also need a comparison principle for for the fractional p-Laplacian (see \cite[Lemma 9]{24}).

\begin{lem}\label{20171096}
Let $\Omega$ be bounded in $R^n$, $p>2$ and $s \in (0,1)$. Assume that
 $u,\,v \in L_{ps}$. If
		\begin{equation}
		\begin{cases}
		(-\Delta)^s_p u\leq 	(-\Delta)^s_p v ,&\ \  x\in \Omega,\\
		u\leq v,&\ \ x\in \Omega^C,
		\end{cases}
		\end{equation}
		then $u\leq v$ in $\Omega$.

\end{lem}

Let's prove Theorem \ref{20171092}.

\textbf{Proof.}
Briefly speaking, the proof consists of two parts. In part one, using the comparison principle we show that $$\|u\|_{L^\infty(\Omega)}<\infty.$$ In part two, we construct an auxiliary function that is Lipschitz continuous near the boundary so as to cover $u(x)$ from above.

 Let $g(x)=\min\{(2-x_n)^s_+,5^s\}$. Then
 \[g(x)=(2-x_n)^s_+-((2-x_n)^s_+-5^s)_+.\]
By \cite[Lemma 3.1]{23}, we know
 $$(-\Delta)^s_p (2-x_n)^s_+=0, \quad x\in B_1. $$ Hence for $ x\in B_1$, we have
\begin{eqnarray*}
&&(-\Delta)^s_p g(x)\\
&=&(-\Delta)^s_p g(x)-(-\Delta)^s_p (2-x_n)^s_+\\
&=&\int_{y_n\leq -3}\frac{[(2-x_n)^s_+-5^s]_+^{p-1}-[(2-x_n)^s_+-(2-y_n)^s_+]^{p-1}}{|x-y|^{n+ps}}dy \\
&=:&I(x).
\end{eqnarray*}
Since $I: B_1(0)\rightarrow \R$ is continuous and positive,  there exists $c>0$ such that
\begin{equation}
(-\Delta)^s_p g(x)\geq c>0 \text{ in }B_1(0).
\end{equation}
Let $\tilde{g}(x)=g(\frac{x}{R})C$ with $R>0, C>0$ sufficiently large so that $\Omega\subset B_R(0)$ and
  $$(-\Delta)^s_p\tilde g(x)= \frac{C^{p-1}}{R^{ps}}[(-\Delta)^s_p g](\frac{x}{R})
  \geq \frac{cC^{p-1}}{R^{ps}}
  \geq \|f\|_{L^\infty(\Omega)} .$$
Then it is obvious that
\begin{equation*}
\left\{\begin{array}{ll}
(-\Delta)^s_p\tilde g(x)\geq (-\Delta)^s_p u(x), & x \in B_R,\\
\tilde{g}(x)\geq u(x), &x \in R^{n}\backslash B_R.
\end{array}
\right.
\end{equation*}
From Lemma \ref{20171096} it follows
$$u(x)\leq \tilde{g}(x) \leq c \text{ in } \Omega.$$
Similarly we can show that
$$-u(x)\leq \tilde{g}(x) \leq c \text{ in } \Omega.$$
This proves that
 $$\|u\|_{L^\infty(\Omega)}\leq c.$$

\medskip

Next we show that $u(x)$ is $C^\nu(\bar{\Omega})$ for $\nu \in (0,s)$. Here
$C^\nu$ denotes the Lipschitz space.
Given $x^o \in \Omega$ and
close to $\partial\Omega$, let $\bar{x^o} \in \partial\Omega$ be such that
$dist(x^o, \partial\Omega)=|x^o\bar{x^o}|$.
We show that
 there  exists a constant $c>0$  such that
\be
|u(x^o)-u(\bar{x^o})|\leq c |x^o-\bar{x^o}|^\nu.
\ee
Without loss of generality, we relocate the origin $O$
so that it is
on the line $x^o\bar{x^o}$ and is outside of $\Omega$ with $|o\bar{x^o}|=1$.
Let $\phi(x)=(|x|^2-1)^\nu_+$. Choose $\xi(x)$ to be a smooth cut-off function
so that $\xi(x)=0$ in $B_1(0)$, $\xi(x)=1$ in  $R^n\backslash B_{1+\epsilon}(0)$ with
the same $\epsilon$ appeared in Lemma \ref{20171098} and $\xi(x)\in [0,1]$ in $R^n$.
Let
$$A(x)=C\phi(x)+\xi(x).$$
Then it is easy to see that $A(x)$ is $C^\nu(\overline{B_1(0)})$.
Without loss of generality,
Let
$$D=B_{1+\epsilon}(0)\backslash B_1(0) \cap \Omega.$$
Given that $\bar{x^o}$ is near $\partial\Omega$, it is reasonable to say that $\bar{x^o} \in D$.

Our goal is to show that
\begin{equation}\label{20171095}
\left\{\begin{array}{ll}
(-\Delta)^s_p A(x)\geq (-\Delta)^s_p u(x),  &x \in D,\\
A(x)\geq u(x), &x \in R^{n}\backslash D.
\end{array}
\right.
\end{equation}
We postpone the proof of (\ref{20171095}) for the moment. Together with Lemma \ref{20171096}, it yields
$$A(x)\geq u(x), \quad x \in D.$$

Since
$$u(x)|_{\partial\Omega}=\xi(x)|_{\partial B_1(0)}=0,$$
and $\xi$ is smooth everywhere, we have
\begin{eqnarray*}
  |u(x^o)- u(\bar{x^o}) |=|u(x^o)|&\leq& |A(x^o)|\\
    &=& |A(x^o)-\xi(\bar{x^o})|\\
     &=& |C (|x^o|^2-1)^\nu_+ + \xi(x^o)-\xi(\bar{x^o})|\\
   &=& |C (|x^o|^2-|\bar{x^o}|^2)^\nu_+ + \xi(x^o)-\xi(\bar{x^o})|\\
   &\leq & C |x^o-\bar{x^o}|^\nu.
\end{eqnarray*}
This implies $u \in C^\nu(\bar{\Omega})$.

What remains is to show (\ref{20171095}). On one hand, it's easy to see that the boundary condition is satisfied because the $A(x)$ controls
$u(x)$ on $R^n \backslash D$ for $C$ sufficiently large. On the other hand, the fractional inequality on
$D$ is valid for $\epsilon$ small because of $\nu<s$ and
\be
(-\Delta)^s_p A(x)\geq C_0(|x|-1)^{\nu(p-1)-ps}, \quad x \in B_{1+\epsilon}(0)\backslash B_1(0).
\ee
To verify this, we use an argument similar to that in the proof of Lemma \ref{20171098}.
 Suppose otherwise,  then there
 exists a sequence $\{x^k\} \in D$ so that
 $|x^k|\rightarrow 1$ and
 \be\label{20171099}
 (-\Delta)^s_p A(x^k)(|x^k|-1)^{\nu(p-1)-ps} \ra 0,  \mbox{ as } \quad  k \ra \infty.
 \ee
Without loss of generality, let $x^k=(0,1+d_k)$. Then
$$d_k=|x^k|-1 \ra 0, \quad \mbox{as } k \ra \infty.$$
By Lemma \ref{c3.1}, we have
\begin{eqnarray*}
 &&(|x^k-1|)^{ps-\nu(p-1)}(-\Delta)_p^s A(x^k)\\
 &=&\frac{d_k^{ps-\nu(p-1)}}{2}C_{n,s,p}\int_{\R^n}
 \frac{[A(x^k)-A(x^k+y)]^{p-1}+[A(x^k)-A(x^k-y)]^{p-1}}{|y|^{n+ps}}dy\\
 &=&C_{n,s,p}d_k^{ps-\nu(p-1)}\bigl(\int_{\R^n}\frac{\Bigl(C(|x^k|^2-1)^\nu_+-C(|x^k+y|^2-1)^\nu_+
 +\xi(x^k)-\xi(x^k+y)\Bigl)^{p-1}}{|y|^{n+ps}}dy\\
 &&\ \ +\int_{\R^n}\frac{\Bigl(C(|x^k|^2-1)^\nu_+-C(|x^k-y|^2-1)^\nu_+
 +\xi(x^k)-\xi(x^k-y)\Bigl)^{p-1}}{|y|^{n+ps}}dy\bigl)\\
&=&\frac{C_{n,s,p}}{2}\int_{\R^n}\bigg(\frac{[C(d_k+2)^\nu_+-C(d_k+2 +2(1+d_k)z_n+d_k|z|^2)^\nu_+
+\frac{\xi(x^k)-\xi(x^k+d_k z)}{d_k^\nu}   ]^{p-1}}{|z|^{n+ps}}\\
 &&\ \ \frac{+\big[C(d_k+2)^\nu_+-C(d_k+2 -2(1+d_k)z_n+d_k|z|^2)_+^\nu
 + \frac{\xi(x^k)-\xi(x^k-d_k z)}{d_k^\nu}
 \big]^{p-1}}{}\bigg)dz
 \end{eqnarray*}
\begin{eqnarray*}
 &=&\frac{C_{n,s,p}}{2}\int_{\R^n}\bigg(\frac{[C(d_k+2)^\nu_+-C(d_k+2 +2(1+d_k)z_n+d_k|z|^2)^\nu_+
+\nabla \xi (\tilde{z})\cdot z d_k^{1-\nu } ]^{p-1}}{|z|^{n+ps}}\\
 &&\ \ \frac{+\big[C(d_k+2)^\nu_+-C(d_k+2 -2(1+d_k)z_n+d_k|z|^2)_+^\nu
 + \nabla \xi (\hat{z})\cdot z d_k^{1-\nu }
 \big]^{p-1}}{}\bigg)dz  \\
 &&\ \ \ \rightarrow\frac{C^{p-1}C_{n,s,p}}{2}
 \int_{\R^n}\frac{(2^\nu-(2+2z_n)_+^\nu)^{p-1}+(2^\nu-(2-2z_n)_+^\nu)^{p-1}}{|z|^{n+ps}}dz\\
 &=&(2^\nu C)^{p-1}\frac{C_{n,s,p}}{2}\int_{\R^n}\frac{(1-(1+z_n)_+^\nu)^{p-1}+(1-(1-z_n)^\nu_+)^{p-1}}{|z|^{n+ps}}dz\\
 &=&(2^\nu C)^{p-1}(-\Delta)_p^s(x_n)^\nu_+|_{x_n=1}\\
 &=&(2^\nu C)^{p-1}C_{\nu,n} >0,
\end{eqnarray*}
 where $\tilde{z}$ is between $x^k$ and $x^k-d_k z$,   $\hat{z}$ is between $x^k$ and $x^k+d_k z$.
This proves (\ref{20171099}) and thus completes the proof of the theorem.


\bigskip

{\em Authors' Addresses and E-mails:}
\medskip

Lingyu Jin

College of Science

South China Agricultural University

Guangdong, Guangzhou, 510640, P.R. China

13822276656@126.com

\smallskip

Yan Li

Department of Mathematics

Baylor University

Waco, Texas, 76706, U.S.

Yan\_Li1@baylor.edu
\end{document}